\documentclass[11pt, reqno]{amsart}

\usepackage{amssymb,latexsym,amsmath,amsfonts,amsthm}
\usepackage{mathrsfs}
\usepackage{enumitem}
\usepackage[usenames]{color}
\usepackage{hyperref}
\usepackage{comment}

\allowdisplaybreaks

\numberwithin{equation}{section}

\theoremstyle{definition}
\newtheorem{definition}{Definition}[section]
\newtheorem{question}[definition]{Question}

\theoremstyle{remark}

\theoremstyle{plain}
\newtheorem{theorem}[definition]{Theorem}

\newtheorem{lemma}[definition]{Lemma}

\newtheorem{corollary}[definition]{Corollary}

\newcommand{\bdy}{\partial}

\newcommand{\opu}{\mathcal{U}}
\newcommand{\nbd}{\mathcal{N}}

\newcommand{\smoo}{\mathcal{C}}

\newcommand{\btl}{\blacktriangleleft}

\newcommand{\Z}{\mathbb{Z}}

\newcommand{\Cn}{\mathbb{C}^n}

\newcommand{\C}{\mathbb{C}} 
\newcommand{\R}{\mathbb{R}}

\begin{document}

\title{Visibility domains that are not pseudoconvex}

\author{Annapurna Banik}
\address{Department of Mathematics, Indian Institute of Science, Bangalore 560012, India}
\email{annapurnab@iisc.ac.in}

\begin{abstract}
The earliest examples of visibility domains, given by Bharali--Zimmer, are pseudoconvex. In fact,
all known  examples of visibility domains are pseudoconvex. We show that there exist non-pseudoconvex
visibility domains. We supplement this proof by a general method to construct a wide range of
non-pseudoconvex, hence non-Kobayashi-complete, visibility domains.
\end{abstract}

\keywords{Visibility, pseudoconvexity, Kobayashi completeness, counterexample, pseudo-arc}
\subjclass[2020]{Primary: 32F45, 32T05; Secondary: 32Q45}

\maketitle

\vspace{-7mm}
\section{Introduction}\label{sec:intro}
The visibility property of domains with respect to the Kobayashi distance has been a subject
of great interest in recent years. If $\Omega\varsubsetneq \Cn$ is a domain such that the Kobayashi
pseudodistance $K_{\Omega}$ is a distance on $\Omega$, then the visibility property, loosely speaking,
is that geodesic lines joining two distinct points in $\bdy{\Omega}$ must bend into $\Omega$ with some mild
geometric control (such as in the Poincar{\'e} disc model of hyperbolic plane). Visibility of
$\Omega$ has been used to deduce many properties of $K_{\Omega}$ and of holomorphic mappings into
$\Omega$ (too
numerous to discuss in this short paper). The visibility property itself has been explored
in several recent papers \cite{bharali-zimmer:2017, bracci-nikolov-thomas:2022,
chandel-maitra-sarkar:2021, bharali-zimmer:2023}.
\smallskip

In this paper, we will focus on $\Omega\varsubsetneq \C^n$ bounded. Our ideas are clearest in this case,
but we will comment on the general case in Section~\ref{sec:lemma_defn}.

\begin{definition} \label{defn:visi-dom}
Let $\Omega \subset \C^n$ be a bounded domain. We say that $\Omega$ is a \emph{visibility domain
with respect to the Kobayashi distance} (or simply a \emph{visibility domain}) if, for any $\lambda \geq 1,
\kappa\geq 0$, each pair of distinct points $\xi, \eta \in \bdy \Omega$, and each pair of
$\overline{\Omega}$-open neighbourhoods $V_{\xi}, V_{\eta}$ of $\xi, \eta$ respectively with $\overline{V_{\xi}}
\cap \overline{V_{\eta}} = \emptyset$, there exists a compact subset $K$ of
$\Omega$ such that the image of any $(\lambda, \kappa)$-almost-geodesic $\sigma:[a,b] \to \Omega$ with
$\sigma(a) \in V_{\xi}, \sigma(b) \in V_{\eta}$ intersects $K$.    
\end{definition}

The visibility property may be viewed as a notion of negative curvature for the metric space
$(\Omega, K_{\Omega})$; it is reminiscent of visibility in negatively curved Riemannian manifolds
in the sense of Eberlein--O'Neill \cite{eberlien-o'neill:1973}. In the latter set-up, the
manifolds\,---\,considered as metric spaces\,---\,are geodesic spaces. In contrast, it is very hard to
tell when the metric space $(\Omega, K_{\Omega})$ is Cauchy-complete\,---\,and hence, when
$(\Omega, K_{\Omega})$ is a geodesic space\,---\,for $n\geq 2$. A
$(\lambda, \kappa)$-almost-geodesic serves as a substitute for a geodesic (since it is known that if
$K_{\Omega}$ is a distance, then for any $z,w \in \Omega$, $z \neq w$, and any $\kappa>0$, there
exists a $(1, \kappa)$-almost-geodesic joining $z$ and $w$ \cite[Proposition~5.3]{bharali-zimmer:2023}).
We refer the reader to Section~\ref{sec:lemma_defn} for the definition of a
$(\lambda, \kappa)$-almost-geodesic.

\smallskip

The earliest examples of visibility domains are the Goldilocks domains, introduced
by Bharali--Zimmer \cite{bharali-zimmer:2017}. They showed that all Goldilocks domains are
pseudoconvex. In fact, all examples of visibility domains in the literature are, or are conjectured
to be, pseudoconvex. This raises the natural question: \emph{are all visibility domains pseudoconvex?}
This question is also of interest as it relates to a deeper issue (which we shall discuss next).
All of this motivates our first theorem:

\begin{theorem} \label{th:non-pscvx_basic}
Let $D \subset \C^n$, $n\geq 2$, be a bounded domain with $\mathcal{C}^2$-smooth boundary and assume
that $D$ is strongly Levi pseudoconvex. Let $A$ be a finite subset of $D$. Then,
$\Omega:= D \setminus A $ is a visibility domain that is not pseudoconvex. In particular,
$\Omega$ is not Kobayashi complete.    
\end{theorem}

The last assertion of Theorem~\ref{th:non-pscvx_basic} has been added for a purpose. (We say
that $\Omega$ is \emph{Kobayashi complete} if $(\Omega, K_{\Omega})$ is Cauchy-complete.) As mentioned above,
the visibility property as given by Definition~\ref{defn:visi-dom} is due, in part, to the fact that it is
very hard to tell when $\Omega$ is Kobayashi complete. Although, in proving results on visibility
domains, $(\lambda, \kappa)$-almost-geodesics are very useful substitutes for geodesics, they lead to proofs
that are sometimes quite technical.
In \cite{bracci-nikolov-thomas:2022} Bracci \emph{et al.} work with bounded convex domains, which are Kobayashi
complete, and with a notion of visibility where $(\lambda, \kappa)$-almost-geodesics are replaced by geodesics
(also see \cite{chandel-maitra-sarkar:2021}).
Also, it follows from \cite[Section~3.2]{chandel-maitra-sarkar:2021} that
if $\Omega$ is Kobayashi complete then, \emph{essentially}, the notion of visibility in
Definition~\ref{defn:visi-dom}, and the notion that \cite{bracci-nikolov-thomas:2022} works with, coincide. For these reasons, experts have
asked
(e.g., at the INdAM conference \emph{New Trends in Holomorphic Dynamics} in Cortona in 2022):
\emph{are all visibility domains Kobayashi complete?} Theorem~\ref{th:non-pscvx_basic} settles this question as well
in the negative. Theorem~\ref{th:non-pscvx_basic} is perhaps the most economical answer to the questions above.
But since the description of the domains in Theorem~\ref{th:non-pscvx_basic} is so specific, one 
may ask whether such examples are an anomaly. Our next result presents a large family of visibility
domains\,---\,which includes the examples in Theorem~\ref{th:non-pscvx_basic}\,---\,that are not Kobayashi complete.

\begin{theorem} \label{th:non-pscvx_Haus-Lip}
Let $D \subset \C^n$, $n\geq 2$, be a bounded visibility domain. Let $A \subset D$ be a compact set with
the following properties:
\begin{enumerate}
  \item $\mathcal{H}^{2n-2}(A)=0$. \label{condn:Haus-measure}
  \smallskip
  
  \item Any Lipschitz map $\phi: I \rightarrow A, \ I$ an interval in $\R$, is constant.
  \label{condn:Lip-path-discnn}
\end{enumerate}
Let $\Omega:= D \setminus A$. Then, $\Omega$ is a visibility domain that is not pseudoconvex. In particular,
$\Omega$ is not Kobayashi complete.
\end{theorem}

Here, for $s \geq 0, \ \mathcal{H}^{s}$ denotes the $s$-dimensional Hausdorff measure. It is very easy to find
totally disconnected sets $A\varsubsetneq D$ having the properties stated in Theorem~\ref{th:non-pscvx_Haus-Lip}.
If these were the only possible cases of $A$, then Theorem~\ref{th:non-pscvx_Haus-Lip} would \emph{essentially} be
a trivial extension of Theorem~\ref{th:non-pscvx_basic} (see the proof of Theorem~\ref{th:non-pscvx_basic}). But
the domains described by Theorem~\ref{th:non-pscvx_Haus-Lip} are a lot more varied because there exist sets $A$
that satisfy the conditions of Theorem~\ref{th:non-pscvx_Haus-Lip} that are connected and contain more than one
point. To see this, we need a definition.
\smallskip

\begin{definition} \label{defn:pseudo-arc} 
A \emph{pseudo-arc} is a compact connected subset $\mathcal{A}\varsubsetneq \R^2$ that contains at least two
points and has the property that no compact connected subset of $\mathcal{A}$ can be expressed as a
union of two non-empty proper subsets that are both compact and connected. 
\end{definition}

Do pseudo-arcs exist? The answer is \textbf{yes}. Given
points $x\neq y\in \R^2$, Moise \cite{moise:1948} gave a construction of a planar set containing the
points $x$ and $y$ and having the properties given in Definition~\ref{defn:pseudo-arc}. In some areas in the
literature, the word ``pseudo-arc'' refers specifically to any planar set homeomorphic to the construction in 
\cite{moise:1948}. (However, we prefer the above definition: it highlights the way in which a pseudo-arc
is \textbf{not an arc,} and Moise's pseudo-arc is covered by this definition.) An equivalent construction of
Moise's pseudo-arc is given in \cite[p.~147]{steen-seebach:1978}. The latter construction is simpler
and establishes the properties given in Definition~\ref{defn:pseudo-arc}.
We can now state a corollary of Theorem~\ref{th:non-pscvx_Haus-Lip}.

\begin{corollary} \label{cor:non-pscvx_simple_Haus-Lip}
Let $D \subset \C^n, n \geq 3$, be a bounded visibility domain. Let $\mathcal{A}\varsubsetneq \R^2$ be
a pseudo-arc and let $A \varsubsetneq D$ be the image of $\mathcal{A}$ under a bi-Lipschitz map
$\psi: \mathcal{A}\to \R^{2n}$. Then,
$\Omega:= D \setminus A$ is a visibility domain that is not pseudoconvex. In particular, $\Omega$
is not Kobayashi complete.
\end{corollary}

In both \cite{moise:1948} and \cite{steen-seebach:1978}, certain inductive
constructions of a pseudo-arc $\mathcal{A}$ are given. In both constructions, $\mathcal{A}$ turns out to be the
intersection of a countable nested family of compact sets that satisfy certain conditions and have non-empty interior.
Due to the last property, it is quite difficult to determine whether $\mathcal{H}^{2n-2}(A)=0$\,---\,where $A$ is
as in Corollary~\ref{cor:non-pscvx_simple_Haus-Lip}\,---\,when $n=2$. For the reader's convenience, we
reproduce the construction of $\mathcal{A}$ from \cite{steen-seebach:1978} in Section~\ref{sec:lemma_defn}. Also, it is
not clear how $\mathcal{H}^{2n-2}(A)$ varies if we consider two inequivalent pseudo-arcs. This is the reason we assume
$n\geq 3$ in Corollary~\ref{cor:non-pscvx_simple_Haus-Lip}; the above uncertainties and ambiguity are avoided,
and $\mathcal{H}^{2n-2}(A)$ is always $0$.

\medskip

\section{Essential lemmas and definitions}\label{sec:lemma_defn}
Most of the definitions in this section do not assume the domain $\Omega \subset
\C^n$ to be bounded. In fact, the theorems in this paper can be established for domains that
are not necessarily bounded (with suitable caveats when $\Omega\varsubsetneq \Cn$ has no
bounded realisation). The lemmas below have substitutes for the unbounded case\,---\,e.g., the
substitute for Lemma~\ref{l:a-g_Lipschitz} is given in \cite{bharali-zimmer:2023}. However,
since we wish to emphasise the underlying \textbf{ideas}, free of the technicalities needed to deal
with the unbounded case, our theorems are stated for bounded domains.
\smallskip

Let $\Omega \subset \C^n$ be a domain. We say that $\Omega$ is \emph{Kobayashi hyperbolic} if $K_{\Omega}$ is an
actual distance. Let $k_{\Omega}: \Omega \times \C^n \to [0, \infty)$ denote the Kobayashi pseudometric. 

\begin{definition} \label{defn:almost-geodesic}
Let $\Omega \subset \C^n$ be a domain and let $I \subset \R$ be an interval. For $\lambda \geq 1$ and
$\kappa \geq 0$, a curve $\sigma: I \to \Omega$ is called a \emph{$(\lambda, \kappa)$-almost-geodesic}
if
\begin{itemize}
  \item ${\lambda}^{-1}|t-s| - \kappa \leq K_{\Omega}(\sigma(s), \sigma(t)) \leq \lambda |t-s| + \kappa
    \quad \forall s,t \in I$, and
   \smallskip
  \item $\sigma$ is absolutely continuous (whereby $\sigma'(t)$ exists for almost every $t \in I$) and
  $k_{\Omega}(\sigma(t); \sigma'(t)) \leq \lambda$ for almost every $t \in I$.
\end{itemize}
    
\end{definition}

We need the following technical lemmas that will play crucial roles in the proof
of Theorem~\ref{th:non-pscvx_Haus-Lip}.

\begin{lemma}[Bharali--Zimmer, {\cite[Proposition~4.3]{bharali-zimmer:2017}}]\label{l:a-g_Lipschitz}
Let $\Omega \subset \C^n$ be a bounded domain. For any $\lambda \geq 1$ there exists a constant
$C=C(\lambda) >0$ such that any $(\lambda, \kappa)$-almost-geodesic $\sigma: I \rightarrow \Omega$
is $C$-Lipschitz (with respect to the Euclidean distance).
\end{lemma}

\begin{lemma} \label{l:non-const_limit_map}
Let $\Omega \subset \C^n$ be a bounded domain. Let $\{ \sigma_{\nu} \}$ be a sequence of
$C$-Lipschitz paths $\sigma_{\nu}:[a_{\nu}, b_{\nu}] \rightarrow \Omega$, for some $C>0$, such that:
\begin{itemize}
  \item There exist constants $a, b \in \R, \ a < 0 < b$, such that $ a_{\nu} \downarrow a$ and
  $b_{\nu} \uparrow b$ as $\nu \to \infty$.
  
\smallskip
  
  \item $\{ \sigma_{\nu} \}$ converges locally uniformly on $(a,b)$ to a continuous map
  $\sigma:(a,b)\rightarrow \overline{\Omega}$.

\smallskip

  \item There exist $\xi', \eta' \in \overline{\Omega}, \ \xi' \neq \eta'$, such that
  $\lim_{\nu \to \infty} \sigma_{\nu}(a_{\nu}) = \xi'$ and $\lim_{\nu \to \infty} \sigma_{\nu}(b_{\nu})=\eta'$.
\end{itemize}
Then, $\sigma$ is non-constant.
\end{lemma}

The proof of the above lemma follows using standard classical analysis; we shall thus omit it. The
essential details of this proof can be found under Case~1 in the proof of \cite[Theorem~1.5]{bharali-maitra:2021}.

\smallskip

We now give a definition that will be helpful in the proof of Theorem~\ref{th:non-pscvx_basic}.

\begin{definition}[Bharali--Zimmer, \cite{bharali-zimmer:2023}] \label{defn:local_M}
Let $\Omega \subset \C^n$ be a Kobayashi hyperbolic domain. Given a subset $\mathcal{S} \subset 
\overline{\Omega}$, we define the function $r \mapsto  M_{\Omega,\,\mathcal{S}}(r) $, $r>0$ as
$$
  M_{\Omega,\,\mathcal{S}}(r) := \sup \bigg\{ \frac{1}{k_{\Omega}(z;v)}: z \in \mathcal{S} \cap \Omega, \, 
  \delta_{\Omega}(z) \leq r, \, \|v\|=1 \bigg\}.
$$
\end{definition}

We end this section with the description given in \cite{steen-seebach:1978} of Moise's pseudo-arc.
What follows is only a \textbf{description;} the reader is referred to \cite{moise:1948, steen-seebach:1978}
for a proof that the set described has the properties mentioned in Definition~\ref{defn:pseudo-arc}. In what
follows, a \emph{chain} $\mathscr{D}$ in $\R^2$ will refer to a finite collection $\{D_j\}_{1\leq j\leq m}$ of
open sets of $\R^2$, called the \emph{links} of $\mathscr{D}$, such that $D_j\cap D_k \neq \emptyset$ if and
only $|j-k|\leq 1$.
Fix two points $x, y \in \R^2$. A pseudo-arc $\mathcal{A} \varsubsetneq \R^2$ containing $x$ and $y$ will
be determined by the sequence of chains $\{\mathscr{D}_{\nu}\}$ constructed inductively as follows:
\begin{itemize}
  \item The diameter of each link of $\mathscr{D}_{\nu} $ is less than $1/{\nu}$.
  \smallskip
  \item The closure of each link of $\mathscr{D}_{\nu + 1}$ is contained in some link of $\mathscr{D}_{\nu}$.
  \smallskip
  \item Let $D^{\nu +1}_{j}, D^{\nu +1}_{k} $ be two links of $\mathscr{D}_{\nu + 1}$, $j < k$. If
  there exist links $D^{\nu}_{\alpha}, D^{\nu}_{\beta}$ of $\mathscr{D}_{\nu}$ with $|\alpha-\beta|> 2$
  such that $ D^{\nu +1}_{j} \subset D^{\nu}_{\alpha} $ and $ D^{\nu +1}_{k} \subset D^{\nu}_{\beta} $,
  then there exist links $D^{\nu +1}_{s}, D^{\nu +1}_{t}$ of $\mathscr{D}_{\nu + 1}$ with $j<s<t<k$
  such that $D^{\nu +1}_{s}$ is contained in some link of $\mathscr{D}_{\nu}$ that is adjacent
  to $D^{\nu}_{\alpha}$ and  $D^{\nu +1}_{t}$ is contained in some link of $\mathscr{D}_{\nu}$
  that is adjacent to $D^{\nu}_{\beta}$.
  
  \smallskip
  \item  For each chain $\mathscr{D}_{\nu}$, its first link contains $x$ and its last link contains $y$.
\end{itemize}
Let $\langle \mathscr{D}_{\nu} \rangle := \bigcup \{D^{\nu}_j: D^{\nu}_{j} \in \mathscr{D}_{\nu} \}$.
Then, $\mathcal{A}:= \bigcap_{\nu} \overline{\langle \mathscr{D}_{\nu} \rangle} $ is a pseudo-arc containing
$x$ and $y$.

\section{The proof of Theorem~\ref{th:non-pscvx_Haus-Lip}} \label{sec:long-proof_non-pscvx}
Given the assumption
\eqref{condn:Haus-measure}, it follows from \cite[Theorem~14.4.5]{rudin:2008} that
$D \setminus A$ is connected. Therefore, by the Hartogs extension phenomenon, the domain
$\Omega=D \setminus A$ is non-pseudoconvex.
Since Kobayashi complete domains are necessarily pseudoconvex
\cite[Theorem~3.4]{kobayashi:2005}, it follows that $\Omega$ is not Kobayashi complete.
\smallskip

We will need the following important observation.

\medskip
\noindent{{\textbf{Claim 1.}} Given $\lambda \geq 1, \kappa  \geq 0$, if $\sigma: [a,b] \rightarrow \Omega$
is a $(\lambda, \kappa)$-almost-geodesic with respect to $K_{\Omega}$, then it is a
$(\lambda, \kappa)$-almost-geodesic with respect to $K_D$ as well.}
\smallskip

\noindent{\emph{Proof of claim:} Given the assumption~(\ref{condn:Haus-measure}}) on $A$,
it follows from \cite[Theorem~3.4.2]{jarnicki-pflug:2013} that
\begin{align}
  K_{\Omega} = K_{D \setminus A} = K_D|_{{(D \setminus A)} \times {(D \setminus A)}} =
  K_D|_{\Omega \times \Omega}. \label{eqn:Kob-dist_equality}
\end{align} 

Using \eqref{eqn:Kob-dist_equality} and the metric-decreasing property for the Kobayashi
metric of the inclusion map $\Omega \hookrightarrow D$ the claim follows immediately from the
definition of a $(\lambda, \kappa)$-almost-geodesic. \hfill $\btl$

\smallskip

We now prove the result by contradiction. Assume that $\Omega$ is not a visibility domain.
Then, there is a pair of distinct points $\xi', \eta' \in \bdy \Omega$, constants $\lambda \geq
1, \kappa \geq 0 $, a pair of sequences $\{z_{\nu}\}, \{w_{\nu}\}$ in $\Omega$ converging to
$\xi', \eta'$ respectively, and a sequence $\{\sigma_{\nu}:[a_{\nu}, b_{\nu}] \rightarrow \Omega \}$
of $(\lambda, \kappa)$-almost-geodesics (with respect to $K_{\Omega}$) joining $z_{\nu}$ and $w_{\nu}$
such that
\begin{equation} \label{eqn:a-g_lim_goes_to_bdry}
  \max_{t \in [a_{\nu}, b_{\nu}] } \delta_{\Omega}(\sigma_{\nu}(t)) \to 0 \ \text{as} \ 
  \nu \to \infty.
\end{equation}

Note that, by our assumptions, $A$ has empty interior. Hence, $\bdy \Omega = \bdy D \cup A$.
So, the following three cases may arise. In each case, we will arrive at a contradiction.

\medskip

\noindent{\emph{Case 1.} $\xi', \eta' \in \bdy D$.}

\smallskip

\noindent{In this case, the fact that $D$ is a visibility domain will play an essential role.}

\smallskip

By Claim~1, each $\sigma_{\nu}$ is a $(\lambda, \kappa)$-almost-geodesic with respect to $K_D$.
Since $D$ is a visibility domain, there is a compact $K$ (that does not depend on $\nu$),
$K \subset D$, such that
\begin{equation} \label{eqn:a-g_intersects_cpt}
  \sigma_{\nu}([a_{\nu}, b_{\nu}]) \cap K \neq \emptyset \ \ \text{for each} \ \nu.
\end{equation}

Since $K \cup A \subset D$ is compact, we can find an open set $\nbd$ such that
$K \cup A \varsubsetneq \nbd \Subset D$. Let us now fix an open neighbourhood $\opu$ of $ \bdy D$
defined as follows:
\[
  \opu := \bigcup\nolimits_{\xi \in \bdy D} B^n(\xi, \epsilon), \quad \text{where}
  \ 0 < \epsilon < \frac{1}{10}\text{dist} \big( \bdy D, \overline{\nbd} \, \big). 
\]
As $ \lim _{\nu \to \infty} 
 \sigma_{\nu}(a_{\nu}) = \xi' \in \bdy D \subset \opu$, there exists $N \in \Z_+$
such that $\sigma_{\nu}(a_{\nu}) \in \opu $ for all $\nu \geq N$.
\smallskip

Consider the auxiliary compact $K' := D \setminus ( \opu \cup \nbd) \subset \Omega$
(note: $K'$ does not depend on $\nu$).

\medskip

\noindent{\textbf{Claim 2.} $\sigma_{\nu}([a_{\nu}, b_{\nu}]) \cap K' \neq \emptyset $ for all
$\nu \geq N$.}

\smallskip

\noindent{\emph{Proof of claim:} Assume that the claim is not true. Then, there is some integer
$\nu' \geq N $ such that $\sigma_{\nu'}([a_{\nu'}, b_{\nu'}]) \cap K' = \emptyset$. This gives
$\sigma_{\nu'}([a_{\nu'}, b_{\nu'}]) \subset \opu \cup \nbd$. Now, $\sigma_{\nu'}(a_{\nu'}) \in 
\sigma_{\nu'}([a_{\nu'}, b_{\nu'}]) \cap \opu $. By \eqref{eqn:a-g_intersects_cpt},
there is a number $t_{\nu'} \in (a_{\nu'}, b_{\nu'}]$ such that $\sigma_{\nu'}(t_{\nu'}) \in 
\sigma_{\nu'}([a_{\nu'}, b_{\nu'}]) \cap K \subset \sigma_{\nu'}([a_{\nu'}, b_{\nu'}]) \cap \nbd .$
\smallskip

But, by our choice of $\opu$, $ \, \overline{\opu} \cap \overline{\nbd} = \emptyset$.
Hence, we conclude from the above that $\sigma_{\nu'}([a_{\nu'}, b_{\nu'}])$ is
disconnected\,---\,a contradiction. Hence the claim. \hfill $\btl$
}
\smallskip

Clearly, Claim~2 contradicts \eqref{eqn:a-g_lim_goes_to_bdry}.

\medskip
\noindent{\emph{Case 2.} $\xi' \in \bdy D, \eta' \in A$.}
\smallskip

\noindent{The investigation of this case is straightforward; it relies on the continuity of
$\sigma_{\nu}$'s.}
\smallskip

Let us consider two open $\epsilon$-neighbourhoods $\opu$ of $\bdy D$ and $\mathcal{V}$ of $A$, where
$0 < \epsilon < 10^{-1}\text{dist} (\bdy D, A)$.
Since $\xi' \in \opu$ and $\eta' \in \mathcal{V}$, by our assumption on $\{\sigma_{\nu}\}$ there is an integer
$N$ such that $\sigma_{\nu}(a_{\nu}) \in \opu $ and $\sigma_{\nu}(b_{\nu}) \in \mathcal{V} $ for all
$\nu \geq N$. Define the compact $K:= D \setminus (\opu \cup \mathcal{V})$; $K\varsubsetneq \Omega$.
Clearly, $\overline{\opu} \cap \overline{\mathcal{V}} = \emptyset$. Thus, as $\sigma_{\nu}([a_{\nu}, b_{\nu}])$
is connected, each $\sigma_{\nu}, \, \nu \geq N$, must intersect $K$. 
This contradicts~\eqref{eqn:a-g_lim_goes_to_bdry}.

\medskip
\noindent{\emph{Case 3.} $\xi', \eta' \in A$.}
\smallskip

\noindent{This case is somewhat involved. By Lemma~\ref{l:a-g_Lipschitz}, there is a constant
$C=C(\lambda)>0$ such that each $\sigma_{\nu}$ is a $C$-Lipschitz map. Then, by the Arzel\`a--Ascoli
theorem, passing to an appropriate subsequence and relabelling if needed, we have:}
\begin{itemize}
  \item $a_{\nu} \downarrow a, \ b_{\nu} \uparrow b$ as $\nu \to \infty$, where $a, b$ are
  in the extended real line.

\smallskip

  \item $\{ \sigma_{\nu} \}$ converges locally uniformly on $(a,b)$ to a continuous map
  $\sigma:(a,b)\rightarrow \overline{\Omega}$.

\smallskip

  \item $ \sigma_{\nu} (a_{\nu})  \to \xi', \  \sigma_{\nu} (b_{\nu}) \to \eta' $
  as $\nu \to \infty$.
\end{itemize}
We now establish the following two claims that will contradict each other.

\medskip

\noindent{\textbf{Claim~3.} The map $\sigma:(a,b)\rightarrow \overline{\Omega}$ is non-constant.}

\smallskip

\noindent{\emph{Proof of claim:} First, we show that $a, b \in \R$. Since each
$\sigma_{\nu}:[a_{\nu}, b_{\nu}] \rightarrow \Omega $ is a $(\lambda, \kappa)$-almost-geodesic,
in view of \eqref{eqn:Kob-dist_equality}, we have }
\begin{equation*}
  \lambda^{-1}|b_{\nu}-a_{\nu}| - \kappa \leq 
  K_{\Omega}(\sigma_{\nu}(a_{\nu}), \sigma_{\nu}(b_{\nu}) ) =
  K_{D}(\sigma_{\nu}(a_{\nu}), \sigma_{\nu}(b_{\nu}) ) \ \ \text{for each } \nu.
\end{equation*}
As $\xi', \eta' \in A$, by
the third bullet point above, \eqref{eqn:a-g_lim_goes_to_bdry}, and using the fact that each
$\sigma_{\nu}$ is continuous, we deduce that there exists $N \in \Z_+$ such that
$\bigcup_{\nu \geq N} \sigma_{\nu}([a_{\nu}, b_{\nu}]) \Subset D$. Hence, $D$ being Kobayashi
hyperbolic, the set $\{ K_{D}(\sigma_{\nu}(s), \sigma_{\nu}(t) ): s,t \in [a_{\nu}, b_{\nu}], \
\nu \in \Z_+ \}$ is bounded above.
Then, the last inequality gives 
$$
  \sup_{\nu}|b_{\nu} - a_{\nu}| < \infty \Rightarrow a,b \in \R.
$$

Also, each $\sigma_{\nu}$ being $C$-Lipschitz, $0 < \|\xi'-\eta' \| \leq C |b-a|$. Hence,
$a \neq b$. Without loss of generality, we may assume $a<0<b$. The claim now follows
from Lemma~\ref{l:non-const_limit_map}. \hfill $\btl$

\medskip

\noindent{\textbf{Claim 4.} The map $\sigma:(a,b) \rightarrow \overline{\Omega}$ is constant.}

\smallskip

\noindent{\emph{Proof of claim:} The second bullet-point above implies that $\sigma_{\nu} \to \sigma$
pointwise. Therefore, $\sigma$ is also $C$-Lipschitz. If, for some $t \in (a,b)$, $\sigma(t) \notin A$,
then it is easy to see, since $\sigma_{\nu}(t) \to \sigma(t)$, that $\delta_{\Omega}(\sigma_{\nu}(t))$
would be bounded away from $0$ for all sufficiently large $\nu$. By \eqref{eqn:a-g_lim_goes_to_bdry},
this is impossible and hence Image$(\sigma) \subset A$.}
Thus, by the assumption \eqref{condn:Lip-path-discnn} on $A$, the Lipschitz map $\sigma:(a,b) \to A$
must be constant.   \hfill $\btl$

\smallskip

Clearly, Claim~3 and Claim~4 contradict each other. Hence, $\Omega$ is a visibility domain, which
completes the proof of the theorem.
\hfill $\qed$

\section{The proofs of Theorem~\ref{th:non-pscvx_basic} and Corollary~\ref{cor:non-pscvx_simple_Haus-Lip}}
\label{sec:short-proof_non-pscvx}
We begin with the short proof of Theorem~\ref{th:non-pscvx_basic}.

\begin{proof}[Proof of Theorem~\ref{th:non-pscvx_basic}]
Note that $\bdy \Omega = \bdy D \cup A $. Let $p \in \bdy \Omega \setminus A = \bdy D.$
We choose a $\C^n$-open neighbourhood $U_p$ of $p$ such that $ U_p \cap \Omega = U_p \cap D$
and $\delta_{\Omega}(z)= \delta_D (z)$ for all $z \in U_p \cap \Omega$. Since
$D$ is a strongly Levi pseudoconvex domain, there is a constant $c>0$ such that
$$
  k_D(z;v) \geq \frac{c \|v\|}{ \sqrt{\delta_D (z)} } \quad \forall z \in D, \ \forall v \in \C^n
$$
(see, for instance, \cite[p.~704]{jarnicki-pflug:2013}). Then, the metric-decreasing
property for the Kobayashi metric of the inclusion map
$\Omega \hookrightarrow D$ and the above choice of $U_p$ give the following
$$
  k_{\Omega}(z;v) \geq   k_D(z;v) \geq \frac{c \|v\|}{ \sqrt{\delta_{\Omega} (z)} }
  \quad \forall z \in U_p \cap \Omega , \ \forall v \in \C^n.
$$
Let $\mathcal{S}_p := U_p\cap \overline{\Omega}$. The above immediately shows that,
for some $\varepsilon >0$, sufficiently small,
\begin{equation} \label{eqn:integrable}
  \int\limits_{0}^{\varepsilon} \frac{1}{r} M_{\Omega,\,\mathcal{S}_p}(r) dr \leq
 \frac{1}{c} \! \int\limits_{0}^{\varepsilon} \! \! \frac{1}{\sqrt{r}} dr 
 < \infty.
\end{equation}

Since $\bdy D$ is $\smoo^2$-smooth and $\bdy \Omega$ is the union of $\bdy D$ with
a finite set $A$, it is easy to prove that $\bdy \Omega$ satisfies an interior-cone
condition in the sense of Bharali--Zimmer \cite[Definition~2.2]{bharali-zimmer:2017}. Then, by
\cite[Lemma~2.3]{bharali-zimmer:2017}, for each $z_0 \in \Omega$, there are constants $C, \alpha >0$
(that depend on $z_0$) such that 
\begin{equation} \label{eqn:Kob-dist_up_bd}
  K_{\Omega}(z, z_0) \leq C + \alpha \log\big(1/{\delta_{\Omega}(z)}\big) \quad \forall z \in \Omega.
\end{equation}
In particular, the above estimate holds for all $z \in U_p \cap \Omega$. Having established
\eqref{eqn:integrable} and \eqref{eqn:Kob-dist_up_bd}, the desired conclusion follows from
\cite[Theorem~1.4]{bharali-zimmer:2023} (also see \cite[Theorem~1.3]{chandel-maitra-sarkar:2021}). To be
more precise: \eqref{eqn:integrable} and \eqref{eqn:Kob-dist_up_bd} imply that each $p\in \bdy D$ is what
is called a ``local Goldilocks point'' in \cite{bharali-zimmer:2023}. Since the set of points in $\bdy \Omega$
that are \textbf{not} local Goldilocks points is totally disconnected, it follows from
\cite[Theorem~1.4]{bharali-zimmer:2023} that $\Omega$ is a visibility domain.
\smallskip

The fact that $\Omega$ is non-pseudoconvex is trivial. Hence, it is not Kobayashi complete.
\end{proof}
We shall now deduce the proof of Corollary~\ref{cor:non-pscvx_simple_Haus-Lip} from
Theorem~\ref{th:non-pscvx_Haus-Lip}.
\begin{proof}[Proof of Corollary~\ref{cor:non-pscvx_simple_Haus-Lip}]
Let $\mathcal{A} \varsubsetneq \R^2$ be a pseudo-arc and let $\psi: \mathcal{A} \to \R^{2n}$ be
a bi-Lipschitz map with $\psi(\mathcal{A})=A \varsubsetneq D$. It is clear that $A$ is a compact
subset of $D$. In order to prove the result, it suffices to show that $A$ satisfies the conditions
\eqref{condn:Haus-measure} and \eqref{condn:Lip-path-discnn} in the statement of
Theorem~\ref{th:non-pscvx_Haus-Lip}.
\smallskip

As $\mathcal{A}$ is planar, its Hausdorff dimension is at most $2$. It follows from the definition of the
Hausdorff dimension that
$$
  \mathcal{H}^s(\mathcal{A}) = 0 \quad \forall s > 2.
$$
Since $n \geq 3$, $\mathcal{H}^{2n-2}(\mathcal{A}) =0$. Since $\psi$ is Lipschitz, by \cite[Proposition~3.1]{falconer:2014}, we have
$$
  \mathcal{H}^{2n-2}(A) = \mathcal{H}^{2n-2}(\psi(\mathcal{A})) \leq C \mathcal{H}^{2n-2}(\mathcal{A}) = 0 \quad
  (\text{for some constant } C>0).
$$
Let $\phi:I \to A$ be a Lipschitz map, where $I \subset \R$ is an interval. Then,
${\psi}^{-1} \circ \phi : I \to \mathcal{A}$ is a continuous path. Assume $\phi$ is non-constant.
Then ${\psi}^{-1}\circ \phi(I)$ contains a compact connected subset of $\mathcal{A}$ that can be
expressed as a union of two non-empty proper subsets that are both compact and connected. In view of
the definition of a pseudo-arc, we have a contradiction. Hence, $\phi$ must be constant.
Thus, both conditions \eqref{condn:Haus-measure} and \eqref{condn:Lip-path-discnn} in the statement of
Theorem~\ref{th:non-pscvx_Haus-Lip} are established for $A$; hence the proof.
\end{proof}

\section{Some questions} \label{sec:questions}
While the results above deal definitively with the questions in Section~\ref{sec:intro},
$\bdy \Omega$, where $\Omega$ is as in the results above, is highly irregular in various senses.
This raises the following:

\begin{question} \label{ques:visi_pscvx_reg-bdry}
Let $\Omega \varsubsetneq \C^n$, $n \geq 2$, be a bounded visibility domain. Suppose the boundary
of $\Omega$ is sufficiently regular\,---\,say, $\bdy \Omega$ is a $\smoo^1$-submanifold. Then, is
$\Omega$ pseudoconvex?
\end{question}

For the reasons discussed in Section~\ref{sec:intro}, a suitable rephrasing of the second question in
Section~\ref{sec:intro} would continue to be of interest. One such rephrasing, which we now discuss, was
suggested by Bracci \cite{bracci:pers-comm2023}.
Since pseudoconvexity and tautness of $\Omega$ are both necessary conditions for $\Omega$ to be Kobayashi
complete, it would be natural\,---\,in view of the theorems above\,---\,to impose one of these conditions
along with visibility to ask whether it follows that $\Omega$ is Kobayashi complete. Since completeness is
already hard to understand, it is reasonable to impose the stronger of the two conditions and ask:

\begin{question} \label{ques:taut_visi_implies_complete}
Let $\Omega \varsubsetneq \C^n$, $n \geq 2$, be a bounded visibility domain. Assume that $\Omega$ is taut.
Then, is $\Omega$ Kobayashi complete?
\end{question}

Furthermore, it is also reasonable to impose the condition of tautness in
Question~\ref{ques:taut_visi_implies_complete} because, for $\bdy \Omega$ sufficiently regular, tautness
and pseudoconvexity coincide; see
\cite{diederich-fornaess:1977, kerzman-rosay:1981, demailly:1987, avelin-hed-persson:2015, chen:2021}.

\section*{Acknowledgements}
\noindent{I would like to thank my thesis advisor, Prof. Gautam Bharali, for introducing to me the questions mentioned
in Section~\ref{sec:intro}. I am also grateful to him for some valuable thoughts and for his feedback on the text of
this paper. I am grateful to Prof. Filippo Bracci for
discussing his interest in the above-mentioned questions. I am also grateful to the anonymous referee of this work
for their help on the exposition above.
This work is supported partly by a scholarship from the
National Board for Higher Mathematics (NBHM) (Ref. No. 0203/16(19)/2018-R\&D-II/10706), and partly by
financial assistance from the Indian Institute of Science.}

\end{document}